\documentclass[12pt]{amsart}
\usepackage{amsbsy,amssymb,amscd,amsfonts,latexsym,amstext,delarray,
amsmath,epsfig}
\input xypic

\newtheorem{thm}{Theorem}[section]  
\newtheorem{prop}[thm]{Proposition}        
\newtheorem{cor}[thm]{Corollary}            
\newtheorem{lem}[thm]{Lemma}                    
 
\newtheorem{defn}{Definition}[section] 
\newtheorem{rem}[thm]{Remark} 
\newtheorem{ex}[thm]{Example}
\numberwithin{equation}{section} 

\def\R{{\mathbb R}}  
\def\H{{\mathbb H}}
\def\Z{{\mathbb Z}} 
\def\Q{{\mathbb Q}} 
\def\C{{\mathbb C}}   
\def\P{{\mathbb P}}
\def\N{{\mathbb N}}
\def\PSL{{\rm PSL}}
\def\PGL{{\rm PGL}}
\def\GL{{\rm GL}}
\def\Red{{\rm Red}}

\begin{document} 

\title[Limiting Modular Symbols and the Lyapunov spectrum]{Limiting
Modular Symbols and the Lyapunov spectrum}  
\author[Matilde Marcolli]{Matilde Marcolli$^{\text{\ddag}}$}
\thanks{$^{\text{\ddag}}$Partially supported by Humboldt Foundation
Sofja Kovalevskaja Award} 
\address{M.~Marcolli: Max--Planck Institut f\"ur Mathematik  \\ 
Vivatsgasse 7 \\ 
Bonn D-53111 \ \ Germany}  
\email{marcolli@mpim-bonn.mpg.de} 
\maketitle 

\section{Introduction}

This paper consists of variations upon a theme, that of {\it
limiting modular symbols} introduced in \cite{MM}. We show that, on 
level sets of the Lyapunov exponent for the shift map $T$ of the continued
fraction expansion, the limiting modular symbol can be computed as a
Birkhoff average. We show that the limiting modular symbols vanish
almost everywhere on $T$--invariant subsets for which a
corresponding transfer operator has a good spectral theory, thus
improving the weak convergence result proved in \cite{MM}. We also
show that, even when the limiting modular symbol vanishes, it is
possible to construct interesting non--trivial homology classes on
modular curves that are associated to non--closed geodesics. These
classes are related to ``automorphic series'', defined in terms of
successive denominators of continued fraction expansion, and their
integral averages are related to certain Mellin transforms of 
modular forms of weight two considered in \cite{MM}. 
We discuss some variants of the Selberg zeta function that sum over
certain classes of closed geodesics, and their relation to Fredholm
determinants of transfer operators. Finally, we argue that one can use
$T$--invariant subsets to enrich the picture of non--commutative
geometry at the boundary of modular curves presented in \cite{MM}.

\section{Modular symbols and geodesics}

Consider the classical compactification of modular curves obtained by
adding cusp points $G\backslash \P^1(\Q)$. Given two points
$\alpha,\beta \in \H\cup \P^1(\Q)$, a real 
homology class $\{ \alpha,\beta \}_G \in H_1(X_G,\R)$ is defined by
integrating the lifts by $\phi : \H \to X_G$ of differentials of
the first kind on $X_G$ along the geodesic arc connecting $\alpha$ and
$\beta$, 
$$ \int_{\{\alpha,\beta \}_G} \omega :=  \int_{\alpha}^\beta \phi^*
\omega. $$

The modular symbols $\{ \alpha,\beta \}_G$ satisfy the additivity and
invariance properties
$$ \{\alpha ,\beta\}_G + \{\beta ,\gamma\}_G = \{\alpha ,\gamma\}_G,
$$
and 
$$
\{g\alpha ,g\beta\}_G=\{\alpha ,\beta\}_G,
$$
for all $g\in G$.
 
Because of additivity, it is sufficient to consider modular
symbols of the form $\{ 0, \alpha \}$ with $\alpha\in \Q$,
$$ \{ 0,\alpha \}_G = -\sum_{k=1}^n \{ g_k(0), g_k(i\infty) \}_G, $$
where $\alpha$ has continued fraction expansion $\alpha =
[a_0,\ldots,a_n]$, and 
$$ g_k =\left(\begin{array}{cc} p_{k-1}(\alpha) & p_k(\alpha) \\
q_{k-1}(\alpha) & q_k (\alpha) \end{array}\right), $$
with $p_k/q_k$ the successive approximations, and $p_n/q_n =\alpha$.

In \cite{MM} it was argued that one can interpret the whole $\P^1(\R)$
with the action of $G$ as a compactification of $X_G$, 
and a corresponding generalization of modular 
symbols was introduced, using infinite geodesics in the upper half
plane $\H$, which end at an {\em irrational} point in $\P^1(\R)$. 

Let $\gamma_\beta$ be an infinite geodesic in $\H$ with one end at
$i\infty$ and the other end at $\beta \in \R \smallsetminus \Q$. Let
$x_0 \in \gamma_\beta$ be a fixed base point, $\tau$ be the geodesic arc
length, and $y(\tau)$ be the point along $\gamma_\beta$ at a distance
$\tau$ from $x_0$, towards the end $\beta$. Let $\{ x,y(\tau) \}_G$
denote the homology class in $X_G$ determined by the image of the
geodesic arc $\langle x, y(\tau) \rangle$ in $\H$. 

The {\it limiting modular symbol} is defined as 
\begin{equation} \label{limitmod} 
\{\{* ,\beta\}\}_G :=\lim \frac{1}{\tau}\,\{ x,y(\tau) \}_G\in
H_1 (X_G,\bold{R}), \end{equation}
whenever such limit exists. The limit \eqref{limitmod} is independent
of the choice of $x_0$ as well as of the choice of the geodesic in
$\H$ ending at $\beta$, as discussed in \cite{MM}. 

\subsection{Shift} We shall study the limiting modular symbols using
properties of the dynamical system given by a generalization of the
shift of the continued fraction expansion. For a modular curve
$$ X_G = \PGL(2,\Z) \backslash (\H \times \P), $$
where $\P = \PGL(2,\Z) /G$, we consider the shift map
$$ T: [0,1] \times \P \to [0,1] \times \P $$
\begin{equation} \label{shift} T(\beta,t) = \left( \frac{1}{\beta} - \left[
\frac{1}{\beta} \right], 
\left(\begin{array}{cc} -[1/\beta] & 1 \\ 1 & 0 \end{array}\right)\cdot t
\right). \end{equation}

\subsection{Lyapunov spectrum}

Recall that the {\em Lyapunov exponent} of the map $T: [0,1] \to [0,1]$
is defined as
\begin{equation} \label{LyapExp} \lambda(\beta) := \lim_{n\to \infty}
\frac{1}{n} \log | (T^n)^\prime (\beta) | = \lim_{n\to \infty}
\frac{1}{n} \log \prod_{k=0}^{n-1} | T^\prime (T^k \beta)
|. \end{equation} 
The function $\lambda(\beta)$ is $T$--invariant. Moreover, in the case of
the classical continued fraction shift $T\beta= 1/\beta -[1/\beta]$ on
$[0,1]$, the Lyapunov exponent is given by
\begin{equation} \label{LyapT} \lambda(\beta) = 2 \lim_{n\to \infty}
\frac{1}{n} \log q_n(\beta), \end{equation}
with $q_n(\beta)$ the successive denominators of the continued
fraction expansion. 
In particular, the Khintchine--L\'evy theorem \cite{Levy} shows that,
for almost all $\beta$'s, the limit \eqref{LyapT} is equal to $2C$,
with  
\begin{equation} \label{C} C=\frac{\pi^2}{12 \log 2}. \end{equation}

The {\em Lyapunov spectrum} is introduced (cf.~\cite{PoWei}) by
decomposing the unit interval in level sets of the Lyapunov exponent
$\lambda(\beta)$ of \eqref{LyapExp}. Let
$L_c = \{ \beta \in [0,1] \, | \lambda(\beta)= c \in \R \}$.
These sets provide a $T$--invariant decomposition of the unit interval,
$$ [0,1] = \bigcup_{c \in \R} L_c \cup \{ \beta \in [0,1] \, |
\lambda(\beta) \text{ does not exist} \}. $$
These level sets are uncountable dense $T$--invariant subsets of
$[0,1]$, of varying Hausdorff dimension \cite{PoWei}. The Lyapunov
spectrum measures how the Hausdorff dimension varies, as a function
$h(c) = \dim_H (L_c)$.

We introduce a function $\varphi: \P \to H_1(\overline{X_G},
cusps, \R)$ of the form 
\begin{equation} \varphi(s) = \{ g(0), g(i\infty) \}_G , \label{varphi}
\end{equation} 
where $g\in \PSL(2,\Z)$ is a representative of the coset $s\in
\P$. We can easily generalize the result \S 2.3 of \cite{MM}, where the
following is proved for the level set $L_{C}$, with $C$ as in
\eqref{C}.  

\begin{thm}
For a fixed $c\in \R$, and for all $\beta \in L_c$,
the limiting modular symbol \eqref{limitmod} is computed by the limit
\begin{equation}\label{limitLyap}
\lim_{n\to \infty} \frac{1}{c n} \sum_{k=1}^n \varphi \circ T^k (t_0),
\end{equation}
where $T: [0,1]\times \P \to [0,1]\times \P$ is the shift operator
\eqref{shift}, and $t_0 \in \P$ is the base point. 
\label{thmLyap}
\end{thm}

\proof  Without loss of generality, we can consider the geodesic
$\gamma_\beta$ in $\H$ with one end at $i\infty$ and the other at
$\beta$. Following the argument given in \S 2.3 of \cite{MM}, we
estimate the geodesic distance $\tau \sim -\log \Im y + O(1)$, as $y(\tau)
\to \beta$. (Here $\Im y$ denotes the imaginary part.) Moreover, if
$y_n$ is the intersection of $\gamma_\beta$ 
and the geodesic with ends at $p_{n-1}(\beta)/q_{n-1}(\beta)$ and
$p_n(\beta)/q_n(\beta)$, then we can estimate
$$ \frac{1}{2 q_n q_{n+1}} < \Im y_n < \frac{1}{2 q_n q_{n-1}}. $$
Moreover, notice that the matrix $g_k^{-1}(\beta)$, with
$$ g_k(\beta)= \left( \begin{array}{cc} p_{k-1}(\beta) & p_k(\beta) \\ 
q_{k-1}(\beta) & q_k(\beta) \end{array}\right), $$
acts on points $(\beta,t)\in [0,1]\times \P$ as the $k$--th power of
the shift operator $T$ of \eqref{shift}. Thus, for $\varphi(t_0)= \{
0,i\infty \}_G$, we obtain
$$ \varphi(T^k t_0) = \{ g_k^{-1}(\beta)\, (0),g_k^{-1}(\beta)\,
(i\infty) \}_G = \left\{ 
\frac{p_{k-1}(\beta)}{q_{k-1}(\beta)},\frac{p_k(\beta)}{q_k(\beta)}
\right\}_G. $$
Finally, we can replace the path $\langle x_0, y_n \rangle$ with the
union of arcs 
$$ \langle x_0, y_0 \rangle \cup \langle y_0, p_0/q_0 \rangle \cup
\bigcup_{k=1}^n \langle p_{k-1}/q_{k-1}, p_{k}/q_{k} \rangle \cup 
\langle p_n/q_n, y_n \rangle $$
representing the same homology class in $H_1(\overline{X}_G,\Z)$. 

\endproof

\subsection{Exceptional set} There is a set of $\beta \in [0,1]$ of
Hausdorff dimension equal to one, where the Lyapunov exponent does not
exist (Theorem 3 of \cite{PoWei}). On this exceptional set  
the limiting modular symbols cannot be expressed as the limit
\eqref{limitLyap} and the corresponding geodesic spends
increasingly long times winding around different closed geodesics. 

\subsection{Periodic continued fractions} A special case of limiting
modular symbols is when the endpoint $\beta$ is a quadratic
irrationality. In this case, $\beta$ has eventually periodic continued
fraction expansion, hence it is a fixed point of some hyperbolic
element in $G$.

Recall that primitive closed geodesics in $X_G$ are parameterized by periodic
points of the shift operator $T$ of \eqref{shift},
$T^\ell (\beta,t) = (\beta,t)$.

\begin{lem}
Consider a geodesic $\gamma_\beta$ in $\H$ with an endpoint at $\beta$
Assume that $\beta$ has eventually periodic continued fraction
expansion. Then the limiting modular symbol is given by 
\begin{equation} \label{period} \{ \{ * , \beta \}\}_G = 
\frac{  \sum_{k=1}^\ell \{
g_k^{-1}(\beta)\cdot g(0), 
g_k^{-1}(\beta)\cdot g(i\infty) \}_G  }{\lambda(\beta) \ell}.
\end{equation} 
\end{lem}

\proof
Consider a geodesic $\gamma_\beta$ in $\H$ with an endpoint at
$\beta$. Without loss of generality we may assume the other endpoint
is at $i\infty$. There is a lift $\gamma_\beta \times \{ t \}$ to
$\H\times \P$, such that $T^\ell (\beta,t) = (\beta,t)$, for a 
minimal non-negative integer $\ell$.

For a quadratic irrationality $\beta$ the limit 
$\lambda(\beta)=2 \lim_{n\to \infty}\log (q_n(\beta))/n$ exists and
belongs to the interval $[2\log((1+\sqrt 5)/2),\infty)$, cf.~ \S 4 of
\cite{PoWei}. 
Thus, we can apply Theorem \ref{thmLyap} and obtain \eqref{period},
where 
$$ \{ g_k^{-1}(\beta)\cdot g(0),
g_k^{-1}(\beta)\cdot g(i\infty) \}_G $$
is the homology class in $\overline{X}_G$ of the geodesic
$$ \pi\left( \langle \frac{p_{k-1}(\beta)}{q_{k-1}(\beta)},
\frac{p_k(\beta)}{q_k(\beta)} \rangle \times \{ g_k^{-1}(\beta) t \}
\right). $$

\endproof

In \cite{MM} a different but equivalent expression for the limiting
modular symbol for quadratic irrationalities is obtained, of the form 
\begin{equation} \{\{ * ,\alpha^+_g \}\}_G = \frac{ \{ 0, g(0) \}_G
}{\lambda(g)}, \label{MSperiodic} \end{equation}
where $\{ 0, g(0) \}_G$ is the homology class of the image in $X_G$ of
the geodesic arc $\langle x_0, g x_0 \rangle$, and $\alpha^\pm_g$ are
the pair of attractive and repelling fixed points of an hyperbolic
element $g\in G$, with $0<\Lambda_g^{\pm}<1$ the respective
eigenvalues, and $\lambda(g) = \log \Lambda^-_g$ is the length of the
closed geodesic in $X_G$ obtained from the geodesic in $\H$ with ends at
$\alpha^\pm_g$.

\section{Transfer and Gauss--Kuzmin operators}

We now treat the case of more general $T$--invariant subsets of
$[0,1]\times \P$. It is often possible to recast the dynamical
properties of a map like our shift \eqref{shift} in terms of the
functional analytic properties of certain transfer operators.

\begin{defn} Let $E\subset [0,1]\times \P$ be a $T$--invariant
subset. The Ruelle (or Perron--Frobenius) transfer operator is defined
as 
$$ ({\mathcal R}_h f)(\beta,t) = \sum_{(\alpha,u) \in
T^{-1}(\beta,t)} \exp(h(\alpha,u)) f(\alpha,u). $$
We only consider the case where the function $h$ is of the form 
$h(\beta,t)=-s \log | T'(\beta) |$, depending on a parameter $s$.
\end{defn}

For a $T$--invariant subset $E\subset [0,1]\times \P$ of Hausdorff
dimension $\delta_E$, we also define a generalized Gauss--Kuzmin
operator, which is the adjoint of composition with $T$ under the
$L^2(E,{\mathcal H}^{\delta_E})$ inner product, with ${\mathcal
H}^{\delta_E}$ the corresponding Hausdorff measure, \cite{Sch}.

\begin{defn} Let $E\subset [0,1]\times \P$ be a $T$--invariant
subset of Hausdorff dimension $\delta_E$. The Gauss--Kuzmin operator
is defined by
\begin{equation} \label{adjoint}
\int_E ({\mathcal L} f)\cdot h \, d{\mathcal H}^{\delta_E} =
\int_E f\cdot  (h\circ T) \, d{\mathcal H}^{\delta_E}, \end{equation}
for all $f,h$ in $L^2(E,d{\mathcal H}^{\delta_E})$.
\end{defn}

Let ${\mathcal N}_E \subset \N$ be the set
\begin{equation} \label{NE} 
{\mathcal N}_E := \left\{ k: \, \, \left(\left[ \frac{1}{k+1},
\frac{1}{k} \right] \times \P \right) \cap E \neq \emptyset
\right\}. \end{equation} 

\begin{rem} {\em
The operator ${\mathcal L}$ defined by \eqref{adjoint}
is given by 
\begin{equation} \label{L1E2}
({\mathcal L} f)(\beta,t) = \sum_{k\in {\mathcal N}_E}
\frac{1}{(\beta+k)^{2\delta_E}} f \left( \frac{1}{\beta+k},
\left(\begin{array}{cc} 0 & 1\\ 1 & k \end{array} \right) \cdot t
\right) , \end{equation}
\label{lemL1E} }
\end{rem}

In this case, we can consider a one--parameter family of deformed
Gauss--Kuzmin operators

\begin{equation} \label{LsE}
({\mathcal L}_{\sigma,E} f)(x,t)  = \sum_{k\in {\mathcal N}_E}
\frac{1}{(x+k)^{\sigma}} f \left( \frac{1}{x+k},
\left(\begin{array}{cc} 0 & 1\\ 1 & k \end{array} \right) \cdot t
\right)  \end{equation} 
With this notation, the operator ${\mathcal L}$ of \eqref{L1E2} is
${\mathcal L}={\mathcal L}_{2\delta_E,E}$.

\begin{rem} {\em If the invariant set is of the form $E=B\times \P$
with $B$ defined by $B = \{ \beta \in [0,1] : a_i \in {\mathcal N} \}$,
for $\beta=[a_1,\ldots, a_n,\ldots]$, and ${\mathcal N}\subset \N$ a
given subset, then we have ${\mathcal R}_h = {\mathcal L}_{2s,E}$,
for $h=-s \log | T'|$. }
\end{rem}

In general, these two operators differ. 
The Ruelle transfer operator tends to carry more 
information on the dynamics of $T$, while the Gauss--Kuzmin 
operator tends to have better functional analytic 
properties. The last Remark describes the optimal case.

\subsection{Spectral theory}

We are especially interested in the cases where the Gauss--Kuzmin
operator and has a good spectral theory.

\begin{defn}
We say that a $T$--invariant set $E\subset [0,1]\times \P$ has good
spectral theory (GST) if the following conditions are satisfied:
\begin{enumerate}

\item There is a real Banach space  ${\mathbb V}$ with the property
that, for each $f\in {\mathbb V}$, the restriction $f|_E$ lies in
$L^2(E,d{\mathcal H}^{\delta_E})$, and the restrictions satisfy
$$ \overline{\{ f|_E \, \, | f\in {\mathbb V} \} }=L^2(E,d{\mathcal
H}^{\delta_E}). $$

\item The operator ${\mathcal L}_{\sigma,E}$ acting on  
${\mathbb V}$ is compact. 

\item There exists a cone $K\subset {\mathbb
V}$ of functions positive at points of $E$, and an
element $u$ in the interior of $K$, such that ${\mathcal 
L}_{\sigma,E}$ is $u$--positive, namely for all non--trivial $f\in K$
there exists $k>0$ and real $a,b >0$ such that
$$ au \leq  {\mathcal L}_{\sigma,E}^k f \leq bu, $$
where the order is defined by $f \leq g$ iff $g-f \in K$.

\item For $\sigma =2\delta_E$ the spectral radius 
of ${\mathcal L}_{2\delta_E,E}$ is one,
$$ \rho( {\mathcal L}_{2\delta_E,E} ) =1. $$

\end{enumerate}
\end{defn}

The techniques developed by Mayer \cite{Mayer} then imply that 
the following properties are satisfied.

\begin{prop}
Let $E$ be a $T$--invariant set with GST. Then For all $\sigma \in I$,
$\lambda_\sigma = \rho({\mathcal L}_{\sigma,E})$
is a simple eigenvalue of ${\mathcal L}_{\sigma,E}$, and the unique
(normalized) eigenfunction satisfying
${\mathcal L}_{\sigma,E} h_\sigma = \lambda_\sigma  h_\sigma$
is in the cone $K$. There exists a unique $T$--invariant measure
on $E$, which is absolutely continuous with respect to the Hausdorff
measure $d{\mathcal H}^{\delta_E}$, with density the normalized
eigenfunction $h_{2\delta_E}$. 

Let ${\mathcal L}_{\sigma,E}^*$ be the adjoint operator acting on
the dual Banach space. There is a unique eigenfunctional $\ell_\sigma$
satisfying ${\mathcal L}_{\sigma,E}^* \ell_\sigma = \lambda_\sigma
\ell_\sigma$. For any $f\in {\mathbb V}$, we have
\begin{equation} \label{convFunctional} \lim_{n\to \infty}
\lambda_\sigma^{-n} {\mathcal L}_{\sigma,E}^n f 
= \ell_\sigma(f) h_\sigma, \end{equation}
with the eigenfunctional $\ell_{2\delta_E}$ given by
$f \mapsto \int_E f\,  d{\mathcal H}^{\delta_E}$.
The iterates ${\mathcal L}_{2\delta_E,E}^k f$, for $f\in
{\mathbb V}$, converge to the invariant 
density $h_{2\delta_E}$, at a rate of the order of $O(q^n)$, where $q$
is the 
spectral margin: $|\lambda| < q \lambda_\sigma$ for all points
$\lambda \neq \lambda_\sigma$ in the spectrum of ${\mathcal
L}_{\sigma,E}$. Thus, for any $f \in {\mathbb V}$ we have
\begin{equation} \label{limAverage}
\lim_{n\to \infty} \frac{1}{n} \sum_{k=1}^n ({\mathcal L}^k f)(\beta,t) =
h_{2\delta_E} (\beta,t) \cdot \int_E f \, d{\mathcal H}^{\delta_E}.
\end{equation}
\label{spectral}
\end{prop}

In particular, given any function $F$ in 
$L^2(E,d{\mathcal H}^{\delta_E})$,  we have
\begin{equation} \label{limIntAverage} 
\int_E \frac{1}{n} \sum_{k=1}^n F(T^k(\beta,t)) \, f(\beta,t) \, d{\mathcal
H}^{\delta_E}(\beta,t) \to \left( \int_E F\, h_{2\delta_E} \,
d{\mathcal H}^{\delta_E}\right) \cdot \left( \int_E f \,
d{\mathcal H}^{\delta_E} \right), 
\end{equation}
for and any test function $f \in
L^2(E,d{\mathcal H}^{\delta_E})$.

\subsection{Generalized Gauss problem}

Let $E=B\times \P$ be a $T$--invariant subset of $[0,1]$ with GST.

Let $t_0 \in \P$ be a base point, and let $m_n(\beta,t)$ be defined as
\begin{equation}\label{mnEt}
m_n(x,t):= {\mathcal H}^{\delta_B}\left( \left\{
\alpha \in B \left| x_n(\alpha)\leq x \, \text{ and } \,
g_n^{-1}(\alpha)\cdot t_0 =t \right. \right\}
\right). \end{equation}

\begin{lem}
The measures \eqref{mnEt} satisfy the recursive relation
$$ m_{n+1}(x,t) = \sum_{k=1}^\infty \left( m_n\left( \frac{1}{k},
\left( \begin{array}{cc} 0 & 1 \\ 1 & k \end{array} \right) \cdot t
\right) 
- m_n \left( \frac{1}{x+k}, \left( \begin{array}{cc} 0 & 1 \\ 1 & k
\end{array} \right) \cdot t \right) \right), $$ 
and the densities satisfy
$$ m_{n+1}'(x,t) = ({\mathcal L}_{2\delta_E,E} m_n') (x,t). $$

The measures $m_n(x,t)$ converge to the unique $T$--invariant measure
$$ m(x,t):= \int_E h_{2\delta_E}(x,t) \, d{\mathcal
H}^{\delta_E}(x,t), $$ 
at a rate $O(q^n)$ for some $0< q <1$.
In this case the density satisfies
$$ h_{2\delta_E}(x,t) = \frac{1}{| \P |} h_{2\delta_B} (x), $$
where $h_{2\delta_B}$ is the top eigenfunction for the Gauss--Kuzmin
operator ${\mathcal L}_{2\delta_B}$ on $B$.
\label{recursionEt}
\end{lem}

The proof of this Lemma is a straightforward generalization of the
proof of Theorem 0.1.2 of \cite{MM}. 

\section{Variation operator}

We introduce another operator, which gives the variation of the
Gauss--Kuzmin operator along the 1--parameter family, ${\mathcal A}_{\sigma,
E} := \frac{d}{d\sigma} {\mathcal L}_{\sigma,E}$. This can be written
in the form
\begin{equation} \label{AsE}
({\mathcal A}_{\sigma,E} f)(\beta,t)= \sum_{k\in {\mathcal N}_E}
\frac{\log(\beta+k)}{(\beta+k)^\sigma} f \left( 
\frac{1}{\beta+k}, \left( \begin{array}{cc} 0 & 1\\ 1 & k
\end{array} \right)\cdot t \right),
\end{equation}
with ${\mathcal N}_E \subset \N$ defined as in \eqref{NE}.

We have the following result.

\begin{lem}
Let $E$ be a $T$--invariant set with GST, and let $\lambda_\sigma$
be the top eigenvalue of the Gauss--Kuzmin operator ${\mathcal
L}_{\sigma,E}$ acting on ${\mathbb V}$, and let $h_\sigma$ be the
corresponding unique (normalized) eigenfunction
\begin{equation}\label{eigenEq}
{\mathcal L}_{\sigma,E} h_\sigma = \lambda_\sigma h_\sigma.
\end{equation}
Then we have
\begin{equation}\label{averageA}
\int_E ({\mathcal A}_{2\delta_E,E} h_{2\delta_E}) \, d{\mathcal H}^{\delta_E}
= \lambda^\prime_\sigma |_{\sigma =2\delta_E},
\end{equation}
with the operator ${\mathcal A}_{\sigma,E}$ as in \eqref{AsE}, and
with $\lambda^\prime_\sigma =\frac{d}{d\sigma} \, \lambda_\sigma$.

\end{lem}

\proof
We differentiate the eigenvalue equation \eqref{eigenEq} with respect
to $\sigma$, and evaluate at $\sigma = 2\delta_E$. This gives
\begin{equation} \label{prime} {\mathcal A}_{2\delta_E,E}
h_{2\delta_E} + {\mathcal 
L}_{2\delta_E,E} h^\prime_{2\delta_E} = \lambda^\prime_{2\delta_E}
h_{2\delta_E} + h^\prime_{2\delta_E}, \end{equation}
where we used $\lambda_{2\delta_E}=1$, and we set 
$$ h^\prime_{2\delta_E} =
\frac{d}{d\sigma} h_\sigma |_{\sigma 
=2\delta_E}. $$
Notice that 
$$ \int_E (h^\prime_{2\delta_E} - 
{\mathcal L}_{2\delta_E,E} h^\prime_{2\delta_E}) d{\mathcal
H}^{\delta_E} =0, $$
hence the result follows by integrating both sides of \eqref{prime}. 

\endproof

Now we return to the
computation of limiting modular symbols \eqref{limitmod}. The
following result complements the result of Theorem \ref{thmLyap}.

\begin{thm}
Let $E=B\times \P$ be a $T$--invariant subset of $[0,1]\times \P$ with
GST. Then, for ${\mathcal H}^{\delta_B}$--almost every $\beta 
\in B$, the Lyapunov exponent satisfies
\begin{equation}\label{Lyaplambda} \lambda(\beta) = 
2\lambda^\prime_{2\delta_B}, \end{equation}
hence the limiting modular symbol \eqref{limitmod} is given by
\begin{equation} \label{limitlambda}  \lim_{n \to \infty}
\frac{1}{2\lambda^\prime_{2\delta_B} \, n} 
\sum_{k=1}^n \varphi(T^k t_0), \end{equation}
where $\varphi$ is defined as in \eqref{varphi}.
\label{eigenLyap}   
\end{thm}

\proof 

The identification \eqref{Lyaplambda} can be proved using a strong law
of large numbers, cf.~\cite{Khin}, \cite{RoSz}. We
compute expectations. For $\beta = [a_1,a_2,\ldots,a_{n-1},a_n + x_n]$,
we have
$$ \frac{1}{x_1 \cdots x_n} = q_{n+1} \left(1 +
\frac{x_{n+2}q_n}{q_{n+1}} \right), $$ 
with $q_n$ the successive denominators satisfying $q_{n+1}=a_{n+1}q_n
+ q_{n-1}$. We can write 
$$ \frac{-1}{n+1} \sum_{k=1}^{n+1} \log (x_n) = \frac{1}{n+1}
\log q_{n+1}  + \frac{1}{n+1} \log \left(1 +
\frac{x_{n+2}q_n}{q_{n+1}} \right). $$ 
Thus, we can estimate the asymptotic behavior of $ \log (q_n) /n$ by the
behavior of the average on the left hand side, as $n\to \infty$.
By the convergence of the measures $m_n$, we have
$\int_B \log(x_n) \, d{\mathcal H}^{\delta} = \int_B  \log(x) \,
h_{2\delta}(x) \, (1 + O(q^n)) d{\mathcal 
H}^{\delta}(x)$.
Here $h_{2\delta}(\beta)$ is the normalized eigenfunction of the
Gauss--Kuzmin operator on $B$, and $\delta_E=\delta_B=\delta$
Moreover, we have
$\int_B  \log(x) \, h_{2\delta} d{\mathcal
H}^{\delta} = \int_B ({\mathcal L}_{2\delta,B} F)\, d{\mathcal
H}^{\delta}$, for
$F(x) := \log(x) h_{2\delta}(x)$.
Notice that by \eqref{AsE} we also have
${\mathcal L}_{2\delta,B} F = -{\mathcal A}_{2\delta,B} h_{2\delta}$,
hence we obtain
$$ \int_B \log(x_n) \, d{\mathcal H}^{\delta} = -\int_B  ({\mathcal
A}_{2\delta,B} h_{2\delta}) (1 + O(q^n)) d{\mathcal H}^{\delta} = 
- \lambda^\prime_{2\delta} (1 + O(q^n)). $$
Thus, we have expectation 
$$ - \int_B \sum_{k=1}^{n+1} \log (x_n) d{\mathcal H}^{\delta} = (n+1)
\lambda^\prime_{2\delta} + O(1). $$
The variances can be computed as
$$ D^2(\sum_{k=1}^{n+1} \log (x_n)) := \int_B \left( \sum_{k=1}^{n+1}
\log (x_n) \right)^2 d{\mathcal H}^{\delta}  -
(\lambda^\prime_{2\delta})^2, $$
and this is estimated by evaluating
$\int_B \log (x_k) \log (x_{k+j}) d{\mathcal H}^{\delta} = \int_B
\log(x) \log(y) dm(x,y)$,
with $dm(x,y)$ the distribution for the measure of $\{ x_k(\beta)\leq
x \text{ and } x_{k+j}(\beta) \leq y \}$. By replacing $x_k$ with a
truncation at some large $N$, $x_k^* =
[a_{k+1},a_{k+2},\ldots,a_{k+N}]$, it is possible to ensure that
$dm(x,y) = dm_k(x) \cdot dm_{k+j}(y) (1 + O(q^j))$.
The argument then follows as in the proof of the classical
Khintchine--Levy theorem \cite{RoSz}. 

\endproof

\section{Vanishing results}

In this section we generalize the weak vanishing result of
\cite{MM}, and improve the convergence in average to convergence
almost everywhere. 

\subsection{Weak convergence}

The result of Proposition \ref{spectral} implies the following vanishing
in the average for limiting modular symbols.

\begin{lem}
Let $E=B\times \P$ be a $T$--invariant set with GST. Then, for all
test functions 
$f \in L^2(E,d{\mathcal 
H}^{\delta_E})$, we have
\begin{equation} \label{weaklimit}
\lim_{\tau\to\infty} \frac{1}{\tau} \int_E \{ x_0, y(\tau) \} \cdot f \,
d{\mathcal H}^{\delta_E} = 0. \end{equation}
\label{vanishing}
\end{lem}  

\proof
The proof is analogous to Lemma 2.3.1 of \cite{MM}.
We write the limiting modular symbol as 
$$ \lim_{n \to \infty} \frac{1}{2 \lambda_{2\delta}' \,n} 
\sum_{k=1}^n \varphi(T^k t). $$ 
Moreover, we have
$$ \int_{B\times \P} \frac{1}{n}\sum_{k=1}^n \varphi\circ T^k \, f
d{\mathcal H}^{\delta_E} = \int_{B\times \P} \frac{1}{n}\sum_{k=1}^n
({\mathcal L}^k f) \varphi \, d{\mathcal H}^{\delta_E}. $$
By \eqref{limIntAverage}, this converges to
$$ \left( \int_{B\times \P} \varphi h_{2\delta_E} \, d{\mathcal
H}^{\delta_E} \right) \left( \int_{B\times \P} f \, d{\mathcal
H}^{\delta_E} \right) $$
where $h_{2\delta_E} =h_{2\delta_B}/|\P |$. This gives
$$ \left( \frac{1}{|\P |}\sum_{s\in \P} \varphi(s) \right)
\left(\int_{B\times \P} f \, d{\mathcal H}^{\delta_E} \right). $$
Arguing as in Lemma 2.3.1 of \cite{MM} we see that $\sum_{s\in \P}
\varphi(s) =0$, hence the result follows.

\endproof

\subsection{Strong convergence}

Now we show that, using a strong law of large numbers it is possible
to improve the weak convergence of Theorem \ref{vanishing} to 
convergence ${\mathcal H}^{\delta_E}$--almost everywhere in $E$. 

\begin{thm}
Let $E=B\times \P$  be a $T$--invariant subset of $[0,1]\times \P$,
with GST. Then, for ${\mathcal
H}^{\delta_E}$--almost every $\beta \in E$ we have
$\{ \{ * , \beta \}\} =0$.
\label{strongvanishing}
\end{thm}

\proof
Here we treat $\varphi_k = \varphi(T^k t_0)$ as random variables, and
prove that we can apply the strong law of large numbers. The argument is
similar to the one used in the proof of Theorem \ref{eigenLyap}.
The result of Theorem \ref{vanishing} implies that the expectation
$E=E(\varphi_k)$ is zero. We evaluate deviations,
$$ D^2 = \int_E \left| \sum_{k=1}^n \varphi(T^k t_0) \right|^2
d{\mathcal H}^{2\delta_E}. $$
Here we write $| \varphi_k |^2$ for the pairing
$| \varphi_k |^2 = \langle \varphi_k, \varphi_k \rangle$, 
where 
$$\langle \varphi(s), \varphi(t) \rangle := \{ g(0),g(i\infty) \}_G
\bullet \{ h(0), h(i\infty) \}_G,$$ with $s=gG$, $t=hG$, and $\bullet$
the intersection product, cf.~ \cite{Mer}.

When writing
$$ \sum_{k=1}^n  \int \left| \varphi_k \right|^2 + \sum_{k=1}^n
\sum_{j=1}^{n-k} \int  \langle \varphi_k, \varphi_{k+j} \rangle, $$
we want to estimate the difference 
\begin{equation} \label{estimD} \int  \langle \varphi_k, \varphi_{k+j}
\rangle - \langle \int 
\varphi_k, \int \varphi_{k+j} \rangle. \end{equation}

Let $P( g_k^{-1}(\beta)\cdot t_0 = t)$ be the probability that
$g_k^{-1}(\beta)\cdot t_0 = t$, that is
$$ P( g_k^{-1}(\beta)\cdot t_0 = t)= \int_B dm_k(x,t) d{\mathcal
H}^{2\delta}(x)= m_k(1,t), $$
for $m_n(x,t)$ defined as in \eqref{mnEt}. We can write \eqref{estimD}
as 
$$ \begin{array}{ll} \sum_{s\in \P} \sum_{t\in \P} \langle \varphi(s),
\varphi(t) \rangle  & ( P(g_k^{-1}(\beta)\cdot t_0 =s \, \text{ and } \,
g_{k+j}^{-1}(\beta)\cdot t_0 = t ) \\[3mm]
&  - P(g_k^{-1}(\beta)\cdot t_0 =s)
\cdot P(g_{k+j}^{-1}(\beta)\cdot t_0 =t) ). \end{array} $$

It is then enough to show that these are sufficiently well approximated
by weakly dependent events, namely that we have
$$ P(g_k^{-1}(\beta)\cdot t_0 =s \, \text{ and } \,
g_{k+j}^{-1}(\beta)\cdot t_0 = t ) = P(g_k^{-1}(\beta)\cdot t_0 =s)
\cdot P(g_{k+j}^{-1}(\beta)\cdot t_0 =t) (1 + O(q^j)). $$ 

By the result of Lemma \ref{recursionEt}, we know that the measures 
$m_n(x,t)$ converge to the limit $T$--invariant measure $m(x,t)$ at a
rate $O(q^n)$. However, the probability distributions of $x_k(\beta)$
and $x_{k+j}(\beta)$ do not satisfy the weak dependence condition
$$ m(x,y) = m_n(x) m_{n+j}(y) (1 + O(q^j)). $$
However, by proceeding as in the proof of Theorem \ref{eigenLyap}, we
consider a truncation at some large 
$N$, $x_k^*=[a_{k+1},a_{k+2},\ldots,a_{k+N}]$, so that the
probabilities 
$$ m_n^*(x,t)= {\mathcal H}^{2\delta}\left( \{ \alpha \in [0,1] \, |
x_k^*(\alpha) \leq x \, \text{and} \, g_k^{-1}(\alpha)\cdot t_0 =t \}
\right) $$
satisfy the weak dependence condition. This gives the desired result
for the probabilities $P(g_k^{-1}(\beta)\cdot t_0 =s \, \text{ and } \,
g_{k+j}^{-1}(\beta)\cdot t_0 = t )$.

\endproof

\section{Examples}

We discuss some examples of invariant sets with GST that are relevant
to the boundary geometry of modular curves.

\subsection{Generalized Gauss--Kuzmin operator} This is the case
analyzed in \cite{MM}. 
In this case the function space ${\mathbb V}$ is the real Banach space
$V({\mathbb D}\times \P)$ of functions holomorphic on each sheet
of ${\mathbb D}\times \P$ and continuous to the boundary, real
at the real points of each sheet, with ${\mathbb D}=\{z\in\C
\, | \, |z-1|<3/2 \}$. The space is endowed with the supremum
norm. The generalized Gauss--Kuzmin operator acts on this
space ${\mathbb V}$, for $\sigma > 1/2$.
An additional hypothesis on the subgroup
$G\subset \PGL(2,\Z)$ is needed, namely the {\em
transitivity condition} that the coset space $\P$
contains no proper invariant subset under the action of the semigroup
$\Red=
\cup_{n\geq 1} \Red_n \subset \GL(2,\Z)$, with
\begin{equation} \label{Red}
\Red_n=\left\{ \,
\left( \begin{array}{cc}
0 & 1\\
1 & a_1
\end{array} \right)\ \dots \
\left( \begin{array}{cc}
0 & 1\\
1 & a_n
\end{array} \right)\,|\,a_1,\ldots , a_n\ge 1;\, a_i\in\Z \right\},
\end{equation}
This ensures that $[0,1]\times \P$ has GST. 
In this case, the $T$--invariant measure on $[0,1]\times \P$ is given by
\begin{equation} \label{GKmeasure} \mu(\beta,t) = \frac{1}{| \P | \log 2}
\log (1+\beta), \end{equation}
where $1/(1+\beta)\log 2 $ is the invariant density of the classical Gauss
problem for the shift $T\beta=1/\beta-[1/\beta]$.
Thus, in this case we obtain an improvement on the weak 
convergence result of \cite{MM}.

\begin{cor}
Assume that $\Red (t)=\P$ for each $t\in\P$. Then the limiting modular
symbol satisfies
$\{\{* ,\beta\}\}_G =\lim_{\tau\to \infty} \,
\frac{1}{\tau}\,\{x_0,y(\tau)\}_G =0$, 
for almost all $\beta\in [0,1]$.
\end{cor}

\subsection{Continued fractions Cantor sets}

There is a family of Cantor sets associated to the continued fraction
expansion of numbers in $[0,1]$, namely the sets $E_N$ given by the
numbers $\alpha\in [0,1]$, $\alpha= [a_1,a_2,\ldots,a_{\ell},
\ldots]$, with all $a_i \leq N$. These sets have Hausdorff dimensions
which tend to one as $N\to \infty$ according to the asymptotic formula
\cite{Hen2}
$$ \delta_N := \dim_H (E_N) = 1- \frac{6}{\pi^2 N} - \frac{72 \log N}{
\pi^4 N^2} + O(1/N^2). $$

As before, let $G$ be a finite index subgroup of $\PGL(2,\Z)$, and
$\P$ be the coset space. There is an action of the shift
operator $T$ restricted to the space $E_N \times \P$, where $E_N$
is one of the Cantor sets as above.

The Gauss--Kuzmin operator is of the form
\begin{equation} \label{truncatedGK} ({\mathcal L}_{\sigma,N} f)(\beta,t) =
\sum_{k=1}^N 
\frac{1}{(\beta+k)^{\sigma}} f\left( 
\frac{1}{\beta+k}, \left( \begin{array}{cc} 0 & 1\\ 1 & k
\end{array} \right)\cdot t \right), \end{equation}
acting on the Banach space $B$ defined in \S 1.1 of \cite{MM}.

In the case $G=\PGL(2,\Z)$ the properties of this operator 
were studied by Hensley in \cite{Hen1}, \cite{Hen2}, where, using a
functional analytic setup similar to Babenko's \cite{Bab}, it is shown
that for $\sigma=2\delta_N$, there is a unique eigenfunction $f_N (\beta)$
of ${\mathcal L}_{2\delta_N,N}$, with eigenvalue $\lambda =1$, which
gives the density of an invariant measure.

Following Hensley, we can consider the truncated Gauss--Kuzmin
operators ${\mathcal L}_{\sigma,N}$ as perturbations of the Gauss--Kuzmin
operator ${\mathcal L}_s$ on $[0,1]$.

\begin{prop}
Under the assumption $\Red (t)=\P$ for each $t\in\P$, and 
for sufficiently large $N\geq N_0$, there exists a
unique eigenfunction $f_N (\beta,t)$ of the operator ${\mathcal
L}_{\delta_N,N}$ of \eqref{truncatedGK}, with eigenvalue $\lambda
=1$. This satisfies 
\begin{equation} \label{eigenfunctionTGK}  f_N (\beta,t) = \frac{1}{|
\P |} f_N(\beta),  \end{equation}
where $f_N(\beta)$ is the unique normalized eigenfunction of the
truncated Gauss--Kuzmin operator of \cite{Hen1}, \cite{Hen2}, for
$\sigma=2\delta_N$, and $\P=\{ 1 \}$.
\label{eigenEN}
\end{prop}

\proof

The assumption $\Red (t)=\P$ for each $t\in\P$ ensures that the
Gauss--Kuzmin operator ${\mathcal L}_\sigma$, for the shift on
$[0,1]\times \P$ has GST for $\sigma > 1/2$.  

For large $N$, we consider the operator ${\mathcal L}_{\sigma,N}$ as a 
perturbation of ${\mathcal L}_\sigma$. We estimate the operator 
norm of $T_{\sigma,N} ={\mathcal
L}_\sigma - {\mathcal L}_{\sigma,N}$. We have
$$ T_{\sigma,N} f(x,t) = \sum_{k=N+1}^\infty k^{-\sigma}
\frac{1}{(x/k+1)^{\sigma}} 
f\left( \frac{1}{x+k}, \left( \matrix 0 & 1\\ 1 & k
\endmatrix \right)\,(t) \right), $$
hence
$$ \| T_{\sigma,N} f \| \leq \sum_{k=N+1}^\infty k^{-\eta} \left\|
\frac{1}{(1 + x/k)^{\sigma}} f\left( \frac{1}{x+k}, \left( \matrix 0 & 1\\ 1
& k \endmatrix \right)\,(t) \right) \right\| $$
$$ \leq \frac{C}{N^{\eta}} \| f \|, $$
for $\sigma$ complex with real part $\Re(\sigma)=\eta >1/2$.

One can then apply the Crandall--Rabinowitz bifurcation Lemma \cite{CrRab}
(cf.~ Lemma 7 of \cite{Hen2}), which guarantees the existence of a $\delta
>0$, such that, for $\| T_{\sigma,N} \| < \delta$, there is a unique
$\rho_{\sigma,N} \in \C$, with $| \rho_{\sigma,N} | < \delta$, for
which $L_\sigma 
-(\lambda_\sigma +\rho_{\sigma,N})I$ is singular. The map
$T_{\sigma,N} \mapsto \rho_{\sigma,N}$ is 
analytic, in the sense specified in \cite{CrRab},
\cite{Hen2}. Moreover, there exists a unique $f_{\sigma,N}$ satisfying 
$$ ({\mathcal L}_\sigma + T_{\sigma,N}) f_{\sigma,N} = (\lambda_\sigma
+\rho_{\sigma,N}) f_{\sigma,N}, $$ 
and the map $T_{\sigma,N} \mapsto f_{\sigma,N}$ is also analytic.
By uniqueness, the eigenfunction for $\sigma=2\delta_N$ is therefore
of the form $f_N (\beta,t) = f_N(\beta)/|\P |$, with 
$f_N(\beta)$ the unique normalized eigenfunction of the
truncated Gauss--Kuzmin operator of Hensley.

\endproof

The eigenfunction $f_N (\beta,t)$ of the truncated Gauss--Kuzmin operator
defines the density of a $T$--invariant measure on $E_N \times
\P$ absolutely continuous with respect to the
$\delta_N$--dimensional Hausdorff measure $d{\mathcal H}^{\delta_N}$. 

We then obtain the following vanishing result.

\begin{lem}
Consider infinite geodesics in the upper half plane with an end
at $i\infty$ and the other end at a point of $E_N$.
For sufficiently large $N\geq N_0$, the limiting modular symbol
satisfies
$$ \{ \{ *, \beta \} \}_G = \lim_n \frac{1}{2 \lambda'_{2\delta_N}\,
n}\sum_{k=1}^n \varphi( g_k(\beta)^{-1} t_0 ) =0. $$
\end{lem}

Notice that the constant $N_0$ provided by the Crandall--Rabinowitz
bifurcation Lemma is independent of $G\subset \PGL(2,\Z)$.  

\section{Intersection numbers}

In the classical theory of the modular symbols of \cite{Man1},
\cite{Mer} (cf.~ also \cite{Ha}, \cite{Maz}) cohomology classes
obtained from cusp forms are evaluated against relative homology
classes given by modular symbols, and the corresponding intersection
numbers are interpreted in terms of special values of $L$--functions
associated to the automorphic forms which determine the cohomology
class. In this and the next section, we seek to extend this viewpoint
to the theory of limiting modular  
symbols. We first recall the setting of \cite{Mer}, which provides a
convenient setting for computing intersection numbers with limiting
modular symbols. We consider finite index subgroups $G\subset
\PSL(2,\Z)$. The arguments given for $\PGL(2,\Z)$ adapt with minor
modifications, cf.~ \cite{LewZa}, \cite{MM}.

Let $I$ and $R$, be the elliptic points on the modular curve
$X_G$, namely the image under the quotient map $\phi: \H \to X_G$ of
the $\PSL(2,\Z)$ orbits of $i$ and $\rho=e^{2\pi i/3}$, respectively. 
With the notation of \cite{Mer}, we set $H_A^B:=
H_1(\overline{X_G}\smallsetminus A, B; \Z)$. The modular symbols $\{
g(0), g(i\infty) \}$, for $gG\in \P$, define classes in $H^{\text{cusps}}$.
For $\sigma$ and $\tau$
the generators of $\PSL(2,\Z)$ with $\sigma^2=1$ and $\tau^3=1$, we
set $\P_I = \langle \sigma \rangle \backslash
\P$ and $\P_R =\langle \tau \rangle \backslash
\P$. There is an isomorphism $\Z^{|\P |} \cong H^R_{\text{cusps} \cup
I}$. Given the exact sequences 
$$ 0 \to H_{\text{cusps}} \stackrel{\iota'}{\to} H^R_{\text{cusps}} 
\stackrel{\pi_R}{\to} \Z^{| \P_R |} \to \Z \to 0 $$
and
$$ 0\to \Z^{| \P_I |} \to H^R_{\text{cusps} \cup I} 
\stackrel{\pi_I}{\to} H^R_{\text{cusps}} \to 0, $$ 
the image $\pi_I(\tilde x)\in H^R_{\text{cusps}}$ of an 
element $\tilde x = \sum_{s\in \P} \lambda_s s$ in $\Z^{|\P|}\cong
H^R_{\text{cusps}\cup I}$ represents an element $x \in
H_{\text{cusps}}$ iff the image $\pi_R(\pi_I(\tilde x)) =0$ in $\Z^{|
\P_R |}$. As proved in \cite{Mer}, for $s=gG \in \P$, the intersection
pairing $\bullet : H^{\text{cusps}} \times H_{\text{cusps}} \to \Z$
gives 
$$ \{ g(0), g(i\infty) \} \bullet x = \lambda_s - \lambda_{\sigma
s}. $$ 
We define the function $\Delta_x : \P \to \R$ by  
\begin{equation} \label{deltax} \Delta_x (s)= \lambda_s -
\lambda_{\sigma s}, \end{equation} 
where $x$ is given as 
above. 

We list some examples of results involving intersection numbers that
easily follow from the general results of the previous section.

\begin{ex} {\em
For a fixed $x\in H_{\text{cusps}}$, represented by a linear
combination $\tilde x = \sum_{s\in \P} \lambda_s s$ in $\Z^{|\P|}$
satisfying $\pi_R(\pi_I(\tilde x)) =0$ in $\Z^{| \P_R |}$, 
the asymptotic intersection number
\begin{equation} \label{limitinters} \{\{ *, \beta \}\} \bullet x :=
\lim_{\tau \to \infty} 
\frac{1}{\tau} \{ x_0, y(\tau) \} \bullet x \end{equation}
is computed by the limit
\begin{equation} \label{limxn} \lim_{n\to\infty}
\frac{1}{\lambda(\beta)\, n} \sum_{k=1}^n 
\Delta_x( g_k^{-1}(\beta) \cdot t_0 ), \end{equation}
where $\lambda(\beta)$ is the Lyapunov exponent. 
Let $E$ be a $T$--invariant set with GST.  
For ${\mathcal H}^{\delta_E}$--almost all $\beta \in E$, 
the sequence \eqref{limxn} converges to  
$$\frac{1}{\lambda |\P|} \sum_{s\in \P} \Delta_x(s) \equiv 0. $$ 
}
\label{limxnlemma}
\end{ex}

In the following example, we consider a cusp form $\Phi$ on $\H$,
obtained as the pullback 
$\Phi=\varphi^*(\omega)/dz$ under the quotient map
$\varphi: \H \to X_G$, and intersection numbers
$\Delta_\omega (s) = \int_{g_s (0)}^{g_s (i\infty)} \Phi(z) dz$, with
$g_s G= s\in \P$. 

\begin{ex} {\em 
Let $\Phi$ be a cusp form on $\H$, and  
let $\beta$ be a number in $[0,1]$ with eventually periodic continued
fraction expansion with period of length $\ell$. Then we have 
\begin{equation} \label{integralxn} \lim_{\tau\to \infty}
\frac{1}{\tau} \int_{x_0}^{y(\tau)} \Phi(z) \, dz =
\frac{1}{\lambda(\beta)\ell} 
\sum_{k=1}^{\ell} \int_{\frac{p_{k-1} (\beta)}{q_{k-1}
(\beta)}}^{\frac{p_{k} (\beta)}{q_{k} (\beta)}} \Phi(z) \, dz =
\frac{1}{\lambda(g)} \int_0^{g(0)} \Phi(z) \, dz. 
\end{equation}
Here $\Phi$ is a cusp form, $\Phi(z) dz = \varphi^*(\omega)$,
with $\omega$ representing the dual of the class $x$. The integral
homology class  $\{ 0, g(0) \}_G$ is determined by the periodic geodesic
with $\beta = \alpha^+_g$, $g\in G$, and with $\lambda(g)= \log
\Lambda^-_g$. }
\end{ex}

The following example follows by arguing as in Theorem 3.5 of \cite{Man1}.

\begin{ex} {\em 
Let $G=\Gamma_0(N)$, and let $\Phi=\varphi^*(\omega) /dz$ be an
eigenfunction for the Hecke operators $T_m$, for $(m,N)=1$, with
rational eigenvalues $c_m$, and let $\beta$ be a quadratic
irrationality with period of length $\ell$ in the continued
fraction. Then we obtain 
$$ \lim_{\tau \to \infty} \int_{x_0}^{y(T)} \Phi(z) dz =
\frac{1}{(\sigma(m)-c_m) \lambda(\beta)\ell} \sum_{d|m} \, \, \sum_{\{ b \mod
d\} } \, \, \int_{\gamma_{m,b,d}(\beta)} \Phi(z) dz, $$
with $\sigma(m)=\sum_{d|m} d$, and for infinitely many $m$ with $(m,N)=1$, 
$$ \gamma_{m,b,d}(\beta):= \sum_{k=1}^\ell \left\{
\frac{p_k(\beta)}{q_k(\beta)} , \frac{m}{d^2} \cdot
\frac{p_k(\beta)}{q_k(\beta)} + \frac{b}{d} 
\right\}_{\Gamma_0(N)} \in H_1(X_{\Gamma_0(N)}, \Z). $$
}
\end{ex}

\begin{ex} {\em 
Let $E$ be a $T$--invariant set with GST. For $g\in \PSL(2,\Z)$,
consider the factor $U(g,z)$ of the form 
$$ U(g,z) = \frac{1}{(cz+d)^2}  \, \, \text{ for }
\, g=\left(\begin{array}{cc} a&b\\ c&d \end{array}\right), $$
and set $\tilde\Phi (z) := \sum_{t\in \P} \Delta_\omega(t)
U(g_t^{-1},z) \Phi(g_t^{-1} z)$. This satisfies 
$$ \lim_{n\to \infty} \frac{1}{n} \sum_{k=1}^n
(\Delta_\omega (g_k^{-1}(\beta) \cdot t_0))^2 \to \frac{1}{|\P |} 
\int_0^{i\infty} \tilde\Phi(z) dz = \frac{1}{|\P |} \sum_{s\in\P}
\Delta_\omega(s)^2, $$ 
for ${\mathcal H}^{\delta_E}$--almost all $\beta\in E$. For instance,
for some $G=\Gamma_0(p)$ with $p$ prime, the value of this limit
can be computed from the tables at the end of \cite{Man1}
}
\end{ex}

\section{Automorphic series on the boundary}

In \cite{MM} it was proved that cusp forms of weight two for
congruence subgroups (or rather their Mellin 
transforms) can be obtained by integrating along the real axis certain
``automorphic series'' defined in terms of continued fractions and
modular symbols. 

In particular, the following identity was proved: for $\Re(t) >0$, 
\begin{equation} \label{idLint} 
\int_0^1 d\beta \sum_{n=0}^{\infty}
\frac{q_{n+1}(\beta)+q_n(\beta)}{q_{n+1}(\beta)^{1+t}}\,
\int_{\{0,\frac{q_n(\beta)}{q_{n+1}(\beta)}\}}\omega  = 
\left[ \frac{\zeta (1+t)}{\zeta (2+t)}
-\frac{L_{\omega}^{(N)}(2+t)}{\zeta^{(N)} (2+t)^2}
\right]\,\int_0^{i\infty}\Phi(z) dz, 
\end{equation} 
where the cusp form $\Phi=\varphi^*(\omega)/dz$ is an eigenform for all
Hecke operators, $L_{\omega}^{(N)}$ is its Mellin transform with
omitted Euler $N$--factor, $\zeta(s)$ is Riemann's zeta, with
corresponding $\zeta^{(N)}$.
Various generalizations of this average are also discussed in
\cite{MM}. The main point in considering such averages is in order to
recast the theory of modular forms in the upper half plane in terms of
some function theory on the ``invisible boundary'' $\P^1(\R)$ of $X_G$. 

The general type of functions considered in \S 2 of \cite{MM} is of
the form
\begin{equation} \label{ell} \ell(f,\beta)=\sum_{k=1}^\infty
f(q_k(\beta),q_{k-1}(\beta)), \end{equation} 
for $f$ a complex valued function defined on pairs of coprime integers
$(q,q')$ with $q\geq q'\geq 1$ and with $f(q,q') =O(q^{-\epsilon})$ for
some $\epsilon >0$. The identity
$$ \int_0^1 d\beta \, \, \ell(f,\beta) = \sum_{\begin{array}{c} q \geq
q'\geq 1 \\ (q,q')=1 \end{array}} \frac{f(q,q')}{q(q+q')}, $$
which is used in \cite{MM} to prove \eqref{idLint} and various
generalizations, is a result of L\'evy, \cite{Levy}.
In our viewpoint, the summing over pairs of
successive denominators is what replaces modularity, when ``pushed to
the boundary''. This correspondence between Dirichlet series
related to modular forms of weight two 
and integral averages of such ``automorphic series'' on the 
``invisible part'' of the boundary of $X_G$ is reminiscent of the
physical principle of holography, or Maldacena correspondence. We
shall return to discuss this relation in a separate work.  

In \cite{MM} only averages $\int_{[0,1]} \ell(f,\beta) d\beta$ were
considered, and their relations to modular forms. Here we
discuss the pointwise behavior, and prove that it is possible to use
such expressions to construct non--trivial homology classes
associated to a geodesic with ``generic'' endpoint in $[0,1]$, and a
function $f$ in this L\'evy class. We shall only discuss the example
of identity \eqref{idLint}, but similar arguments hold for the other
functions considered in \cite{MM}.

\begin{thm}
Consider the function
\begin{equation} \label{fqq} f(q,q')=\frac{q+q'}{q^{1+t}}
\int_{\left\{ 0, \frac{q'}{q} \right\}_{\Gamma_0(N)}}
\omega. \end{equation}
Let $\ell(f,\beta)$ be the corresponding function as in \eqref{ell}. 
We can estimate
$$ \ell(f,\beta) \sim \sum_{n=1}^\infty e^{-(5+2t)n \lambda(\beta)}
\sum_{k=1}^n \Delta_\omega( g_k^{-1}(\beta)\cdot t_0 ). $$
Thus, if $E$ be a $T$--invariant subset with GST, and with
$0< \lambda^\prime_{2\delta_E}< \infty$, the series defining
$\ell(f,\beta)$ converges absolutely, for 
${\mathcal H}^{\delta_E}$--almost all $\beta \in E$. 
\label{convergence1}
\end{thm}

\proof We can estimate
$$ \frac{q_n(\beta)+q_{n+1}(\beta)}{q_{n+1}(\beta)^{1+t}} \sim
e^{-(5+2t)n \lambda(\beta)}, $$
with $\lambda(\beta)=2 \lim_{n\to \infty} \log(q_n)/n$.  
Moreover, by equation (2) of \S 1 of \cite{Man1}, we can write
$$ \{ 0, \frac{b}{a} \} = \sum_{k=-1}^n \{ \frac{b_{k-1}}{a_{k-1}},
\frac{b_{k}}{a_{k}} \} = \sum_{k=-1}^n \{ g_k^{-1}(b/a) (0),
g_k^{-1}(b/a) (i\infty) \}, $$
where $b_k/a_k$ are the successive convergents of the continued
fraction expansion of the rational number $b/a=b_n/a_n$.

Recall that, in terms of continued fractions, we can write the rational
numbers 
$$ \frac{q_n(\beta)}{q_{n+1}(\beta)}=[a_n,a_{n-1},\ldots,a_1], $$
where
$$ \beta=[a_1,a_2,\ldots,a_{n-1},a_n,a_{n+1},\ldots]. $$
Since the successive denominators $q_n$, viewed as polynomials 
$$ q_n(\beta) = Q_n (a_1,\ldots,a_n), $$
satisfy
$$ Q_n (a_1,\ldots,a_n) = Q_n (a_n,\ldots, a_1), $$
we have, for $\xi_n(\beta):=q_n(\beta)/q_{n+1}(\beta)$,
$$ g_k^{-1}(\xi_n(\beta)) =g_k^{-1}(\beta), $$
for $k\leq n$. Hence the estimate on $\ell(f,\beta)$ follows.
Then arguing as in Theorem \ref{strongvanishing} we obtain
$$ \left| \sum_{k=1}^n \Delta_\omega( g_k^{-1}(\beta)\cdot t_0 )
\right| \sim o(\lambda(\beta) n), $$
hence the series converges absolutely.

\endproof

In fact, this proves a stringer result, namely that we can construct
this way non--trivial homology classes depending on the point $\beta$
and on the function $f$. 

\begin{thm}
Let $E$ be a $T$--invariant subset of $[0,1]\times \P_{\Gamma_0(N)}$
with GST, and with $0< \lambda_{2\delta_E}^\prime < \infty$. 
Then, for ${\mathcal H}^{\delta_E}$--almost all $\beta \in E$, and
for $\Re(t)>0$, the limit
\begin{equation} \label{Cbeta} C(f,\beta):=
\sum_{n=1}^\infty \frac{q_{n+1}(\beta)+
q_n(\beta)}{q_{n+1}(\beta)^{1+t}}\,  
\left\{ 0,\frac{q_n(\beta)}{q_{n+1}(\beta)}\right\}_{\Gamma_0(N)}
\end{equation} 
defines a class $C(f,\beta)$ in $H_1(X_{\Gamma_0(N)},\R)$. 
The function $\ell(f,\beta)$, for $f$ as in \eqref{fqq}, can be
written as
$$ \ell(f,\beta) = \int_{C(f,\beta)} \omega. $$   
\label{convergence2}
\end{thm}

\proof For $G=\Gamma_0(N)$, we have
$$ \sum_{n=1}^\infty \frac{q_{n+1}(\beta)+
q_n(\beta)}{q_{n+1}(\beta)^{1+t}}\,  
\left\{ 0,\frac{q_n(\beta)}{q_{n+1}(\beta)}\right\}_{G} \sim
\sum_{n=1}^\infty e^{-(5+2t)n \lambda(\beta)}
\sum_{k=1}^n \left\{ g_k^{-1}(\beta)(0), g_k^{-1}(\beta)(i\infty)
\right\}_G. $$
Then the same argument used in Theorem \ref{convergence1} proves the
result.

\endproof

The average $\int_{[0,1]} d\beta \int_{C(f,\beta)} \omega$ satisfies
\eqref{idLint}. 
It would be interesting to know if averages over other $T$--invariant
subsets, $\int_E d{\mathcal H}^{\delta_E}(\beta) \int_{C(f,\beta)}
\omega$ also carry number theoretic significance.

\section{Selberg zeta}

The Selberg zeta function of a modular curve $X_G$, for $G$ a finite
index subgroup of $\PGL(2,\Z)$, can be recovered
from the Fredholm determinant of the generalized Gauss--Kuzmin
operator, namely 
\begin{equation}\label{selberg} \det (I-{\mathcal L}_{2s}) = Z_G(s),
\end{equation}
for $s\in \C$ with $\Re(s) > 1/2$. The result holds for $G\subset
\PSL(2,\Z)$ with ${\mathcal L}_{2s}$ replaced by ${\mathcal
L}_{2s}^2$, see \cite{ChM}, 
\cite{MM} for the general case, and \cite{LewZa} for the case
$G=\PSL(2,\Z)$ or $\PGL(2,\Z)$.  

We investigate possible generalizations of this identity when only
geodesics with ends on some smaller $T$--invariant subset are
considered. 

Given a $T$--invariant subset $E\subset [0,1]\times \P$, we can define
\begin{equation} \label{zetaEgeod} 
Z_{G,E}(s):= \prod_{\gamma\in {\rm Prim}_E} \prod_{m=0}^\infty  \left(
1-e^{-(s+m) \text{length}(\gamma)} \right), \end{equation}
where ${\rm Prim}_E$ is the set of primitive closed geodesics in $X_G$
that lift to geodesics in $\H \times \P$ with ends at points of
$E$. Notice that, in general, little is known about the properties of
such zeta functions. For instance, whether they still have a
meromorphic continuation, or if there is a trace formula relating
\eqref{zetaEgeod} to the spectrum of the Laplacian on some suitable
space of functions.

Regarding the relation of \eqref{zetaEgeod} to the Gauss--Kuzmin
operator, the following result is a simple generalization of
the results of \cite{MM} and \cite{Mayer}. It applies, for instance,
to Hensley's Cantor sets of continued fractions $E_N\times \P$.

\begin{prop} 
Let $E$ be a $T$--invariant subset of $[0,1]\times \P$ with GST, which
is of the form $E=B\times \P$, with
\begin{equation} \label{BN} B = \{ \beta\in [0,1] : a_i \in {\mathcal
N} \} \end{equation} 
for some ${\mathcal N} \subset \N$. We have a trace formula
\begin{equation} \label{detzeta} \det (I-{\mathcal L}_{2s,E}) =
Z_{G,E}(s), \end{equation} 
where ${\mathcal L}_{2s,E}={\mathcal R}_h$ is the Ruelle transfer
operator, with $\sigma=2s \in \C$, $\Re(s)>1/2$, $h=-s \log |T'|$, 
and the right hand side is defined as in \eqref{zetaEgeod}.
\label{traceform}
\end{prop}

\proof

We denote by $\pi_{\sigma,k}$ the operator
$$ (\pi_{\sigma,k} f)(x,t) = \frac{1}{(x+k)^{\sigma}} f\left(
\frac{1}{x+k}, \left(\begin{array}{cc} 0&1\\1&k \end{array} \right) \cdot t 
\right), $$
for $\sigma \in \C$, $\Re(\sigma)>1$. We can write \eqref{LsE} as the
operator 
\begin{equation} \label{LsE2} 
{\mathcal L}_{\sigma, E} f  = \sum_{k\in {\mathcal N}}
\pi_{\sigma,k} f. \end{equation}

The same argument of \cite{Mayer} shows that this is a nuclear
operator of order zero in the sense of \cite{Groth} on a complex
Banach space of holomorphic functions defined on a domain ${\mathbb
D}\times \P$ containing $E$, for $\Re(\sigma)>1$. Thus, we have
$$ -\log \det (I-{\mathcal L}_{\sigma,E})= \sum_{\ell =1}^\infty
\frac{Tr\, {\mathcal L}_{\sigma,E}^\ell}{\ell} = Tr\left( \sum_{\ell
=1}^\infty \frac{1}{\ell} \left( \sum_{k\in {\mathcal N}}
\pi_{\sigma,k} \right)^\ell \right) = $$
$$ Tr\left( \sum_{g \in \Red_E} \frac{1}{\ell(g)} \pi_{\sigma}(g)
\right) = \sum_{g\in {\rm Hyp}_E} \frac{1}{\kappa(g)} \chi_\sigma(g)
\tau_g, $$
where $\pi_{\sigma}(g) = \prod_{j=1}^{\ell(g)} \pi_{\sigma,a_j}$, and 
$$ \Red_E := \left\{ g \in \Red \left| \, g =\prod_{j=1}^{\ell(g)}
\left(\begin{array}{cc} 0 & 1 \\ 1 & a_j \end{array} \right)
\right. \, a_i \in {\mathcal N} \right\}. $$
Here we denote by ${\rm
Hyp}_E$ a set or representatives of the 
conjugacy classes of hyperbolic matrices in $\GL(2,\Z)$ which contain
reduced representatives in $\Red_E$, and
$$ \chi_{2s}(g) = \frac{N(g)^{-s}}{1-\det(g) N(g)^{-1}}, $$
for
$$ N(g)= \left( \frac{Tr(g) + D(g)^{1/2}}{2} \right)^2, $$
and $D(g)=Tr(g)^2 -4\det(g)$, and $\tau_g$ is defined as
$\tau_g:= \# \{ t\in \P : g\cdot t =t \}$.

The argument now follows the lines of \cite{LewZa}. 
Under the hypotheses on $E$, \eqref{zetaEgeod} is of the form
\begin{equation} \label{zetaE} 
Z_{G,E}(s):= \prod_{g\in {\rm Prim}_E} \prod_{m=0}^\infty \det \left(
1-\det(g)^m N(g)^{-(s+m)} \rho_{\P}(g) \right), \end{equation}
where $\rho_{\P}(g)$ is the action of $g$ on $\P$, and ${\rm Prim}_E$
is the set of 
primitive hyperbolic elements, $g = h^{k(g)}$, $g\in {\rm Hyp}_E$.
We write \eqref{zetaE} as
$$ \sum_{g\in {\rm Prim}_E} \sum_{m=0}^\infty Tr \sum_{k=1}^\infty
\frac{1}{k} \det(g)^{mk} N(g)^{-(s+m)k} \rho_{\P}(g^k)= $$
$$ \sum_{g\in {\rm Prim}_E} \sum_{k=1}^\infty \frac{1}{k}
\frac{N(g)^{-ks} \tau_{g^k}}{1-\det(g^k) N(g)^{-k}}=
\sum_{g\in {\rm Hyp}_E} \frac{1}{\kappa(g)} \chi_\sigma(g)
\tau_g, $$
where $Tr(\rho_{\P}(g))=\tau_g$. The argument underlying the passage
from summing over ${\rm Prim}_E$ to the summing over ${\rm Hyp}_E$ is
exactly as in \cite{LewZa}, by observing that if a conjugacy class in ${\rm
Hyp}$ contains a representative in $\Red_E$, then all the
$\ell(g)/k(g)$ representatives in $\Red$ are also in $\Red_E$. 
Thus, we obtain 
$$ -\log \det (I-{\mathcal L}_{\sigma,E})= -\log Z_{G,E}. $$

\endproof

It is easy to see where this argument breaks down when less strong
assumptions are made on the invariant set $E$. 
For instance, for a $T$--invariant subset $E=B\times \P$, we can
consider the sets $\Red_{E,\ell} = \left\{ g\in \Red_\ell :
[\overline{a_1,\ldots, a_\ell}] \in B \right\}$. Even assuming that
all the corresponding operators $\pi_{\sigma}(g)$ are of trace class,
we no longer have the identification of the term $\sum_{g\in
\Red_{E,\ell}} \pi_{\sigma}(g)$ with $(\sum_n \pi_{\sigma,n})^\ell$,
hence we can no longer identify $Tr ( \sum_{\ell=1}^\infty
\frac{1}{\ell} \sum_{g\in \Red_{E,\ell}} \pi_{\sigma}(g))$ with the
Fredholm determinant of a Gauss--Kuzmin operator.
Many interesting examples of invariant subsets, such as level sets of
the Lyapunov exponent, do not satisfy the assumptions of Proposition
\ref{traceform}. However, in certain cases, it is still possible to
identify \eqref{zetaEgeod} with a Fredholm determinant, albeit for a
different Gauss-Kuzmin operator.

\subsection{Example: minimal Lyapunov exponent}

The minimal value attained by the Lyapunov exponent at periodic points
of the shift $T$ on $[0,1]$ is given by 
$$ c_0 = 2 \log \left( \frac{1 +\sqrt 5}{2} \right). $$
Let $L_{c_0}$ be the level set of the Lyapunov exponent
$L_{c_0} =\{ \beta\in [0,1] : \lambda(\beta)=c_0 \}$.

\begin{prop}
Consider the dense $T$--invariant subset of $L_{c_0}\times \P \subset
[0,1]\times \P$. The zeta function \eqref{zetaEgeod}, with
${\rm Geod}_E$ the set of primitive closed geodesics in $X_G$ which
lift to geodesics in $\H$ with end in $L_{c_0}$ satisfies 
$$ \det (I-{\mathcal L}_{2s,E_1\times\P}) =Z_{G,L_{c_0}\times\P}(s), $$
for $s\in \C$, $\Re(s) >1/2$, where ${\mathcal
L}_{\sigma,E_1\times\P}$ is the transfer operator associated to the
Hensley Cantor set $E_N\times\P$ for $N=1$. 
\end{prop}

\proof

The level set $L_{c_0}$ of the Lyapunov exponent consists of all
$\beta \in [0,1]$ whose continued fraction expansion has a proportion
of $1$'s which tends asymptotically to $100 \%$, cf.~\cite{PoWei}.
Thus, $L_{c_0}$ is an uncountable dense subset of $[0,1]$. The only
eventually periodic points in $L_{c_0}$ are numbers with partial
quotients in the continued fraction expansion eventually equal to
$1$. Thus, we have
$$ -\log Z_{G,L_{c_0}\times \P}(s) = Tr \left( \sum_{\ell=1}^\infty
\frac{1}{\ell} \sum_{g\in \Red_{L_{c_0}\times \P,\ell}}
\pi_{\sigma}(g) \right) = \sum_{\ell=1}^\infty
\frac{1}{\ell} Tr \, {{\mathcal L}_{2s,E_1}}^\ell, $$
for
$$ ({\mathcal L}_{2s,E_1} f)(x,t)= \frac{1}{(x+1)^{2s}} f \left(
\frac{1}{x+1}, \left(\begin{array}{cc} 0&1\\1&1
\end{array}\right)\cdot t \right). $$

\endproof

A case which is especially of interest, in relation to Manin's theory
of real multiplication \cite{Man2}, and which also does not fit
into the assumptions of Proposition \ref{traceform} is the case of
geodesics in $\H$ with ends at quadratic irrationalities in a real
quadratic field $\Q(\sqrt d)$. One can separate out the Selberg zeta
function into contributions of different real quadratic fields and
seek for corresponding Gauss--Kuzmin operators. We shall return to
this in future work.

\section{$C^*$--algebras}

We show how $T$--invariant subsets of $[0,1]\times \P$ can
be used to enrich the picture of non--commutative geometry at the
boundary of modular curves, already illustrated in \cite{MM}. 

Consider a $T$--invariant set $E=B\times \P$, where $B$ is defined by
a condition on the digits \eqref{BN}, for ${\mathcal N}\subset \N$.
Such sets are totally disconnected compact $T$--invariant subsets of
$[0,1]\times \P$. Let $C(E,\Z)$ be the space of continuous functions
from $E$ to the integers. We denote by $C(E,\Z)^T$ the  invariants 
$C(E,\Z)^T= \{ f\in C(E,\Z) \,  | \, f -
f\circ T =0 \}$,  
and by $C(E,\Z)_T$ the co--invariants $C(E,\Z)_T = C(E,\Z)
/B(E,\Z)$, with the set of co--boundaries 
$B(E,\Z):=\{ f - f\circ T \, | \, f\in
C(E,\Z) \}$. 

We recall the following result (see e.g.~\cite{BoHa}).

\begin{prop}
Let $E$ be a totally disconnected
compact Hausdorff space. There is an identification
$K_0(C(E)) \cong C(E,\Z)$,
while $K_1(C(E))=0$. Moreover, there is an exact sequence
$$ 0 \to C(E,\Z)^T \to K_0(C(E)) \stackrel{ I-T_*}{\to} K_0(C(E)) \to 
C(E,\Z)_T \to 0. $$
The set of co--invariants $C(E,\bold{Z})_T$ is a
unital pre--ordered group, 
$(C(E,\Z)_T, C(E,\Z)_T^+, [1])$. 
\label{CET} 
\end{prop}

The pre--order structure is given by specifying the positive cone
$C(E,\Z)_T^+$ and an order unit, namely an element $u\in
C(E,\Z)_T^+$ such that, for all $[f]\in C(E,\Z)_T$ there
exists a $n\in \Z$, $n\geq 1$, such that $nu - [f] \in
C(E,\Z)_T^+$. Here we set:
$$ C(E,\Z)_T^+ :=\{ [f]\in C(E,\Z)_T \, | \exists f_0 \in
C(E,\Z), [f_0]=[f], f_0 \geq 0 \}, $$
with order unit $u =[1]$, the class of the constant function on
$E$.

The pre--ordered groups for $T$--invariant $E$ defined by \eqref{BN}
give a collection of invariants describing some features of the
non--commutative geometry at the boundary of the modular curve
$X_G$. The arithmetic relevance of these invariants remains to be
understood. The structure of these pre--ordered groups can be studied
through traces.

\subsection{Invariant measures and traces}

Recall that a {\it trace} on a unital pre--ordered group 
$(C(E,\Z)_T,C(E,\Z)_T^+,[u])$, where $E$ is a totally disconnected
compact Hausdorff space invariant under $T$,
is a positive homomorphism
$\tau : C(E,\Z)_T \to \R$, where positive means that it 
sends the cone $C(E,\Z)_T^+$ to $\R^{\geq 0}$, which moreover
satisfies $\tau(u)=1$. Normalized Borel $T$--invariant measures on $E$
define traces by integration, 
$$ \tau([f]) = \int_E f d\mu. $$

Thus, if the invariant set $E$ has GST, so that we can prove the
existence of an invariant measure, we have an associated trace on 
$C(E,\Z)_T$ 
$$ f \mapsto \int_E f \, h_{2\delta_E,E} \, d{\mathcal
H}^{\delta_E}. $$

\subsection{Cuntz--Krieger algebras}

Consider a $T$--invariant set $E=B\times \P$, where $B$ is defined by
the condition \eqref{BN}, with ${\mathcal N}\subset \N$ a finite
subset. Then we can further enrich the picture of the non--commutative
geometry at the boundary of the modular curve $X_G$, by considering a
$C^*$--algebra associated to the dynamics of the shift $T$ on $E$, in
the form of a Cuntz--Krieger algebra ${\mathcal O}_{A_E}$ for a matrix
$A_E$, see \cite{Cu}, \cite{CuKri}. 

\begin{ex} {\em
In the case of the Hensley Cantor sets, the shift $T$ on $E_N$ is a
full one--sided shift on $N$ symbols. For $\P=\{ 1 \}$, the associated
$C^*$--algebra is the classical Cuntz algebra ${\mathcal O}_N$. In
general, we consider the Markov partition given by
$$ {\mathcal A}_N := \{ ((i,t),(j,s)) | U_{i,t} \subset T(U_{j,s}) \},
$$ 
where $i,j \in \{ 1, \ldots, N \}$, and $s,t\in \P$, with sets
$U_{i,t} = U_i \times \{ t \}$, where $U_i \subset E_N$ are the clopen
subsets where the local inverses of $T$ are defined,
$$ U_j = \left[\frac{1}{j+1},\frac{1}{j}\right]\cap E_N. $$
The condition $U_{i,t} \subset T(U_{j,s})$ corresponds to
$$ \left(\begin{array}{cc} 0 & 1 \\ 1 & j \end{array} \right) \cdot t
=s . $$
The Markov partition determines a corresponding $Np\times Np$ matrix
$A_N$, for $p=|\P |$, with entries $(A_N)_{is,jt} =1$ if 
$U_{i,t} \subset T(U_{j,s})$ and zero otherwise. We obtain a
Cuntz--Krieger algebra ${\mathcal O}_{A_N}$. }
\end{ex}

There are numerical invariants that can be extracted from the
non--commutative geometry of such Cuntz--Krieger algebras, and our
hope is to show that some of these recover the rich arithmetic
structure of modular curves. An instance where some of this structure
is recovered using $C^*$--algebras is the result of \cite{MM}, where
we showed that from the Pimsner--Voiculescu exact sequence for the
$C^*$--algebra $C(\P^1(\R)\times \P)\rtimes \PSL(2,\Z)$ one can
recover Manin's presentation of the modular complex. 

\medskip

\noindent {\bf Acknowledgement.} Conversations with Yuri Manin 
were a great source of inspiration.

\end{document}